\documentclass[a4paper,12pt,reqno]{amsart}
\usepackage{amssymb}
\usepackage[dvips]{graphicx}
\usepackage{hyperref}
\usepackage[all]{xy}
\usepackage{url}

\input epsf.tex
\newdimen\xsize
\newdimen\oldbaselineskip
\newdimen\oldlineskiplimit
\def\restorelineskip{\baselineskip=\oldbaselineskip%
\lineskiplimit=\oldlineskiplimit}
\def\putm[#1][#2]#3{
\hbox{\vbox to 0pt{\parindent=0pt%
\vskip#2\xsize\hbox to0pt{\hskip#1\xsize $#3$\hss}\vss}}}%
\long\def\Line#1{\hbox to \hsize{#1}}
\def\putt[#1][#2]#3{
\vbox to 0pt{\noindent\hskip#1\xsize\lower#2\xsize%
\vtop{\restorelineskip#3}\vss}}

\makeatletter
\def\xbig[#1]#2{{\hbox{$\m@th\left#2\vbox to#1\xsize{}%
\right.\n@space$}}}
\makeatother
\def\xlar[#1]#2{%
\smash{\mathop{ \hbox to #1\xsize{\leftarrowfill}}\limits^{#2}}}
\def\xrar[#1]#2{%
\smash{\mathop{ \hbox to #1\xsize{\rightarrowfill}}\limits^{#2}}}
\def\xline[#1]{\hbox to #1\xsize{\leaders\hrule\hfill}}

\DeclareFontFamily{U}{rsf}{\skewchar\font'177}%
\DeclareFontShape{U}{rsf}{m}{n}{<-6>rsfs5<6-8>rsfs7<8->rsfs10}{}%
\DeclareFontShape{U}{rsf}{b}{n}{<-6>rsfs5<6-8>rsfs7<8->rsfs10}{}%
\DeclareMathAlphabet\RSFS{U}{rsf}{m}{n}
\SetMathAlphabet\RSFS{bold}{U}{rsf}{b}{n}
\def\sf#1{{\mathsf{#1}}}

\def\slsf{\slshape \sffamily }

\def\msmall#1{\mathchoice{\hbox{\small$\displaystyle {#1}$}}{#1}{#1}{#1}}

\let\xrar=\xrightarrow

\def\cc{{\mathbb C}}

\def\rr{{\mathbb R}}
\def\sss{{\mathbb S}}

\def\pp{{\mathbb P}}

\def\adyn{\sf{1}}

\def\cos{\sf{cos}\,}
\def\const{\sf{const}}

\def\codim{\sf{codim}\,}

\def\dim{\sf{dim}\,}

\def\sym{\sf{Sym}}

\def\grad{\sf{grad}}

\def\i{\sf{i}}

\def\im{\sf{Im}\,}

\def\lim{\mathop{\sf{lim}}}

\def\pr{\sf{pr}}

\def\pr{\sf{pr}}

\def\r{\sf{r}\,}

\def\sin{\sf{sin}\,}
\def\Sing{\sf{Sing}\,}

\def\v{{\mathrm{v}}}

\def\n{{\mathrm{n}}}

\def\sym{\sf{Sym}\,}

\def\eps{\varepsilon}

\def\<{\langle}\let\la=\<
\def\>{\rangle}\let\ra=\>
  
\def\comp{\Subset}
\def\d{\partial}

\def\ddef{\mathrel{{=}\raise0.3pt\hbox{:}}}
\def\deff{\mathrel{\raise0.3pt\hbox{\rm:}{=}}}

\def\fraction#1/#2{\mathchoice{{\msmall{ #1\over#2}}}%
{{ #1\over #2 }}{{#1/#2}}{{#1/#2}}}

\def\le{\leqslant}

\def\emptyset{\varnothing}

\def\longpoints{\leaders\hbox to 0.5em{\hss.\hss}\hfill \hskip0pt}
\def\stateskip{\smallskip}
\def\state#1. {\stateskip\noindent{\bf#1. }} 
\def\statep#1. {\stateskip\noindent{\bf#1 }} 
\def\proof{\state Proof. }

\def\Chi{\raise 2pt\hbox{$\chi$}}

\def\ie{\hskip1pt plus1pt{\sl i.e.\/,\ \hskip1pt plus1pt}}

\def\sli{{\sl i)} } 
\def\slii{{\sl i$\!$i)} } 
\def\sliii{{\sl i$\!$i$\!$i)} }

\def\Chi{\raise 2pt\hbox{$\chi$}}

\let\phI=\phi\let\phi=\varphi\let\varphi=\phI
\def\eps{\varepsilon}

\def\comp{\Subset}
\def\d{\partial}

\def\1{{1\mkern-5mu{\rom l}}}

\def\ge{\geqslant}

\def\fraction#1/#2{\mathchoice{{\msmall{ #1\over#2}}}%
{{ #1\over #2 }}{{#1/#2}}{{#1/#2}}}

\def\le{\leqslant}

\def\emptyset{\varnothing}

\def\qed{\ \ \hfill\hbox to .1pt{}\hfill\hbox to .1pt{}\hfill $\square$\par}

\def\comment#1\endcomment{}

 \belowdisplayskip=\abovedisplayskip

\def\lineeqqno(#1){\hfill\llap{\vbox to 10pt%
{\vss\begin{align} \eqqno(#1)\end{align}\vss}}\vskip1pt}

\advance\headheight 1.2pt

\def\ShowwLLabel#1{}

\def\thechpt{\Roman{chpt}}

\def\newchapt[#1]#2{\newpage%
\refstepcounter{chpt}\setcounter{subsection}{0}%
\setcounter{thm}{0}\setcounter{defi}{0}%
\setcounter{rema}{0}\setcounter{exrc}{0}%
\renewcommand{\thesubsection}{\thechpt.\arabic{subsection}}%
\section*{\begin{center}\huge \bf Chapter \thechpt\\
#2 \end{center}}\label{#1}%
\ \smallskip%
\markboth{Chapter \thechpt}{#2}%
}
\def\newsect[#1]#2{\refstepcounter{section}\setcounter{equation}{0}%
\renewcommand{\thesubsection}{\arabic{section}.\arabic{subsection}}%
\section*{\arabic{section}.
#2}\vspace{-20pt}\label{#1}\vspace{20pt}%
\markboth{Section \arabic{section}}{#2}}

\def\newlect[#1]#2{\refstepcounter{section}%
\renewcommand{\thesubsection}{\arabic{section}.\arabic{subsection}}%
\section*{Lecture \arabic{section}\\
#2}\label{#1}%
\markboth{Lecture \arabic{section}}{#2}}

\def\newprg[#1]#2{\refstepcounter{subsection}%
\subsection*{{\thesubsection.\ #2}} \label{#1}%
}

\def\newappx[#1]#2{%
\refstepcounter{appx}\setcounter{section}{0}%
\renewcommand{\thesubsection}{A\arabic{appx}.\arabic{subsection}}%
\section*{Appendix \arabic{appx}\\ #2}
\label{#1}%
\markboth{Appendix A\arabic{appx}}{#2}
}

\newtheorem{thm}{Theorem}[section]
   \def\newthm#1{\begin{thm}\label{#1}}

\newtheorem{nnthm}{Theorem}
   \def\newthm#1{\begin{nnthm}\label{#1}}

\newtheorem{lem}{Lemma}[section]
   \def\newlemma#1{\begin{lem} \label{#1}}

\newtheorem{prop}{Proposition}[section]
   \def\newprop#1{\begin{prop}\label{#1}}

\newtheorem{nnprop}{Proposition}
   \def\newprop#1{\begin{nnprop}\label{#1}}

\newtheorem{corol}{Corollary}[section]
   \def\newcorol#1{\begin{corol} \label{#1}}

\newtheorem{nncorol}{Corollary}
   \def\newcorol#1{\begin{nncorol} \label{#1}}

\newtheorem{defi}{Definition}[section]
   \def\newdefi#1{\begin{defi} \label{#1}\rm }

\newtheorem{exmp}{Example}[section]
   \def\newexmp#1{\begin{exmp} \label{#1}\rm }

\newtheorem{nnexmp}{Example}
   \def\newexmp#1{\begin{nnexmp} \label{#1}\rm }

\newtheorem{exrc}{Exercise}
   \def\newexrc#1{\begin{exrc} \label{#1}\rm }

\newtheorem{rema}{Remark}[section]
   \def\newrema#1{\begin{rema} \label{#1}\rm }

\newtheorem{nnrema}{Remark}
   \def\newrema#1{\begin{nnrema} \label{#1}\rm }

\def\eqqno(#1){\label{(#1)}}
\def\eqqref(#1){(\ref{(#1)})}

\title{One side continuity of meromorphic mappings\\ between real analytic hypersurfaces}
\author{S. Ivashkovich}
\date{\today}

\address{
Universit\'e de Lille-1, UFR de Math\'ematiques, 59655 Villeneuve
d'Ascq, France} \email{ivachkov@math.univ-lille1.fr}
\address{IAPMM Nat. Acad. Sci. Ukraine
Lviv, Naukova 3b, 79601 Ukraine}


\begin{document}
\begin{abstract}
We prove that a meromorphic mapping, which sends a peace of a real analytic 
strictly pseudoconvex hypersurface in $\cc^2$ to a compact subset of $\cc^N$
which does not contain germs of non-constant complex curves is continuous from 
the concave side of the hypersurface. This implies the analytic continuability 
along CR-paths of germs of holomorphic mappings from real analytic hypersurfaces 
with non-vanishing Levi form to the locally spherical ones in all dimensions.
\end{abstract}

\maketitle

\setcounter{tocdepth}{1}
\tableofcontents

\newsect[INTROD]{Introduction}

\newprg[res]{Statement of the result}
Let $0\in M\subset U$ be a real analytic strictly pseudoconvex hypersurface in a
neighborhood $U$ of zero of $\cc^2$  defined as $M=\{z: Q(z,\bar z) =0\}$, where $Q$ is a 
real analytic function in $U$ with the non-vanishing gradient. In an appropriate local 
coordinates we can suppose that $Q(z,\bar z) = y_2 -|z_1|^2 - r(z_1,\bar z_1, x_2)$, where 
$r$ is of order $\ge 3$ at zero, see \eqqref(cm-nf) below. Set 
\[
U^{+} = \{ z\in U: y_2 > |z_1|^2 + r(z_1,\bar z_1, x_2)\} 
\]
and call this open set the {\slsf concave}  side of $M$.

\smallskip Let furthermore $K\subset U'$ be a compact in a complex manifold $U'$. Our goal 
in this paper is to prove the following

\begin{nnthm}
\label{mer-cont}
Let $M$ and $K$ be as above and suppose, in addition, that $K$ does not contain germs of 
non-constant complex curves. Let $f:U \to U'$ be  a meromorphic mapping  
such that $f|_M[M]\subset K$. Then $f$ is continuous on $\overline{U^+}$.
\end{nnthm}
A meromorphic mapping $f:U\to U'$ between complex manifolds $U$ and $U'$ is defined by its
graph $\Gamma_f$, which is a locally irreducible analytic subset of $U\times U'$ such that
the natural projection $\pi : \Gamma_f\to U$ is proper and generically one to one. If
$\pi : \Gamma_f\to U$ is (proper and) generically $d$ to $\adyn$ then $f$ is called a 
$d$-valued meromorphic
correspondence. In the case of $U'=\pp^N$ a meromorphic map $f$ is defined by a couple 
$f=(f_1,...,f_N)$ where $f_j$ are meromorphic functions, see section \ref{INTER} for more
details. $f|_M[M]$ is by definition the closure of the set $\{f(z):z\in M\setminus I_f\}$,
where $I_f$ denotes the set of points of indeterminacy of $f$. Therefore in our case 
condition $f|_M[M]\subset K$ means that for every $z\in M\setminus I_f$ one has
that $f(z)\in K$. 

\smallskip Notice that in the case $K\comp \cc^N\subset \pp^N$ condition $f|_M[M]\subset K$ 
easily (by maximum principle) implies that $f(\overline{U^+})\subset \hat K$, where $\hat K$ 
is the polynomial hull of $K$, \ie that $f$ is bounded from the concave side $\overline{U^+}$
of $M$. But we claim more: that $f$ is continuous from this side up to $M$.
 
\begin{nnrema} 
\label{minimal} \rm
A good example of a compact without germs of complex curves is a compact real analytic 
hypersurface in $\cc^{n'}$, see \cite{DF}.
\end{nnrema}

\newprg[appl]{Applications}  Let us explain the interest in such a theorem. Recall the following 
result of Pinchuk, see Theorem 6.2  of \cite{P}. Every germ of a holomorphic mapping from a real 
analytic hypersurface $M\subset \cc^n$ to the unit sphere $\sss^{2n'-1}\subset \cc^{n'}$ analytically 
extends along any $CR$-path in $M$. A $CR$-path in $M$ is a path $\gamma : [0,1]\to M$ such that 
$\dot\gamma (t)\in T^c_{\gamma (t)}M$ for every $t\in [0,1]$. The proof of this theorem in \cite{P} 
consist of two steps. First, one proves that $f$ extends {\slsf meromorphically} along $\gamma$, 
see Lemma 6.7 there. Then one proves the holomorphicity of the extended map, see Lemma 6.6 in \cite{P}. 
The proof of both lemmas in \cite{P} crucially uses the assumption that $M'=\sss^{2n'-1}$ and 
does not hold already for a general locally spherical hypersurface on the place of $\sss^{2n'-1}$.

\begin{nnrema} \rm 
It is claimed  in \cite{SV} that a germ of a holomorphic mapping $f:(M,x)\to (M',x')$ of a strictly 
pseudoconvex hypersurface $M\subset \cc^n$ to a compact, real algebraic, strictly pseudoconvex 
hypersurface $M'\subset \cc^{n'}$ analytically extends along any $CR$-path in $M$. Unfortunately 
the proof of \cite{SV} contains a serious gap in Lemma 3.4. Example of \cite{IM} is actually a 
couterexample to this proof. However a careful inspection of \cite{SV} yields a {\slsf meromorphic} 
extension of $f$. We attach an Appendix to our paper in order to help the interested reader to see 
that the proof of \cite{SV} gives the following: {\it in the conditions as above $f$ meromorphically 
extends along any $CR$-path in $M$ starting at $x'$.} See Theorem \ref{SVT} in the Appendix.
\end{nnrema}
So, following the
logic of \cite{P} we are interested in proving that the extended {\slsf meromorphic} map is actually {\slsf holomorphic}, 
\ie we want to have an analog of Lemma 6.6 from  \cite{P} in a possibly more general case. We can do that for 
locally spherical $M'$-s. Namely, combining our theorem with the result of Pinchuk we get the following

\begin{nncorol}
\label{hol-ext1}
Let $(M,x)$ be a germ of a real analytic s.p.c. hypersurface in $\cc^n$, $n\ge 2$ and 
$M'\subset \cc^{n'}$ a compact locally spherical hypersurface. Let $f:(M,x)\to (M',x')$ 
be a germ of a meromorphic mapping. Then $f$ is holomorphic.
\end{nncorol}

Abbreviation s.p.c. stands for strictly pseudoconvex. 
Indeed, let $F:\widetilde{M}'\to \sss^{2n'-1}$ be the {\slsf development} map constructed 
in \cite{BS}. Here $\pi:\widetilde{M}'\to M'$ stands for the universal cover of $M'$. Using our theorem (after
the reduction to the dimension two), we can localize the problem from the concave side of $M$ and then
apply the result of Pinchuk to $F\circ \pi^{-1}\circ f$. This gives us the holomorphicity of $f$.

\smallskip Remark that we do not need algebraicity of $M'$ here. Combining this result with already mentioned
extract from \cite{SV} we get the following
 
\begin{nncorol}
\label{hol-ext2}
Let $(M,x)$ be a germ of a real analytic hypersurface s.p.c. in $\cc^n$, $n\ge 2$ 
and $M'\subset \cc^{n'}$ a compact real algebraic locally spherical hypersurface. Let $f:(M,x)\to (M',x')$ 
be a germ of a holomorphic mapping. Then $f$ analytically extends along every $CR$-path in $M$.
\end{nncorol}

Since the distribution of complex tangents on a strictly pseudoconvex 
hypersurface is contact one can draw a $CR$-path between any two points in $M$, see \cite{G}. Therefore
we obtain one more corollary. 

\begin{nncorol}
\label{hol-ext3}
Suppose that hypersurface $M\subset \cc^2$ is compact real analytic s.p.c. and simply
connected. Let $f:(M,x)\to (M',x')$ be a germ of a non-constant holomorphic map into a compact real algebraic 
s.p.c. hypersurface $M'\subset \cc^{n'}$. Then $f$ extends to a proper holomorphic mapping $F:D\to D'$ between domains 
bounded by $M$ and $M'$ respectively. Moreover $F$ is continuous up to the boundary.
\end{nncorol}

We do not know whether $F$ is holomorphic in a neighborhood of $\bar D$ in the non-spherical case, this is an open question. 
But $F$  {\slsf is meromorphic} in a neighborhood of $\bar D$.

\begin{nnrema} \rm 
{\bf 1)} Algebraicity of $M'$ in Corollaries \ref{hol-ext1}, \ref{hol-ext2} is needed already for 
the meromorphic extension of $f$, see example in \cite{BS}. 

\smallskip\noindent{\bf 2)} When $M'$ is just a compact without the germs of complex curves as in our theorem
one cannot hope to make $f$ continuous also from the convex side of $M$ (this would imply that $f$ is actually 
holomorphic). A counterexample was given in \cite{IM}.

\smallskip\noindent{\bf 3)} The result of Lemma 6.6 from \cite{P} was later reproved in \cite{C} in a 
slightly more general case of $M'$ of the form 
\begin{equation}
\eqqno(model)
M'=\{z'\in \cc^{n'}:\sum_{j=1}^{n'}|z_j|^{2m_j}<1\},
\end{equation}
where each $m_j$ is a positive integer. But let us remark that mapping $\Phi:(z_1,...,z_{n'})\to 
(z_1^{m_1},...,z_{n'}^{m_{n'}})$ sends this $M'$ to $\sss^{2n'-1}$ and $\Phi$ is proper. Therefore 
the result of \cite{C} follows from that of \cite{P}.
\end{nnrema}

\newprg[INTROD.corr]{Case of meromorphic correspondences}

Now suppose that in the conditions of Theorem \ref{mer-cont} our $f$ is a $d$-valued
meromorphic correspondence with values in a complex manifold $U'$. When saying that $f$ on $M$ 
takes its values in $K\subset U'$ we mean that for every $m\in M\setminus I_f$ all values of $f$ at $m$ 
are contained in $K$. By saying that $f$ is continuous on $\overline{U^+}$ we mean that the 
restriction $f|_{\Gamma_f\cap \overline{U^+}}:\Gamma_f\cap \overline{U^+} \to \overline{U^+}$ is 
finite to one everywhere. 

\begin{nncorol}
\label{corr-cont}
Let $M$ and $K$ be as in Theorem \ref{mer-cont}  and let $f:(M,0)\to U'$ be a germ of 
a meromorphic correspondence such that $f|_M[M]\subset K$. Then $f$ is continuous on 
$\overline{U^+}$.
\end{nncorol}

The proof uses the same ingredients as that of Theorem \ref{mer-cont} plus a simple observation 
that $d$-valued meromorphic correspondence from $U$ to $U'$ can be viewed as a meromorphic
mapping from $U$ to $\sym^d(U')$, here $\sym^d(U')$ is the $d$-th symmetric power of $U'$,
and that $\sym^d(K)$ does not contain germs of non-constant complex curves, see section \ref{APL}
for more details.

\smallskip\noindent{\slsf Acknowledgements.} I'm grateful to Rasul Shafikov and Alexander
Sukhov for the stimulating discussions on the subject of this paper.

\newsect[INTER]{Intersections of analytic disks with a real analytic s.p.c. hypersurface}

\newprg[INTER.gen]{Generalities}

Let $(M,0)$ be a germ of a strictly pseudoconvex hypersurface in $\cc^2$ defined in 
some sufficiently small neighborhood $U$ of the origin as the zero set of a real 
analytic function $M\deff\{z\in U: Q (z,\bar z)=0\}$, $Q (0)=0$ and $\grad_0 Q\not=0$.
In an appropriate coordinates $z_1=x_1+iy_1, z_2=x_2+iy_2$ of $\cc^2$ the 
defining function of $M$  can be written as
\begin{equation}
\eqqno(eqm-1)
Q (z, \bar z) = y_2-|z_1|^2 -r(z_1,\bar z_1,x_2).
\end{equation}
This represents our hypersurface $M$ in the form
\begin{equation} 
\eqqno(eqm-2)
M = \{(z_1,z_2): y_2 =  |z_1|^2 + r(z_1,\bar z_1,x_2)\},
\end{equation}
where the rest $r$ can be supposed to be put into the Chern-Moser normal form
\ie has the form 
\begin{equation}
\eqqno(cm-nf)
r(z_1,\bar z_1, x_2) = \sum_{k,l\ge 2}r_{k\bar l}(x_2)z_1^k\bar z_1^l.
\end{equation}
Here $r_{kl}(x_2)$ are real analytic functions on $x_2$ with $r_{kl}=\bar r_{lk}$, 
see \cite{CM}. In fact $r_{k\bar l}$ in the Chern-Moser normal form satisfy some
additional properties, but we shall not need them. Set $U^{\pm} \deff \{z\in U: Q(z,
\bar z)\gtrless 0\}$ and call $U^{\pm}$ the concave/convex side of $M$.

\smallskip Meromorphic mapping $f$ from a complex manifold $U$ with values in a complex
manifold (or  complex space) $U'$ is given by an irreducible and locally irreducible 
analytic subset $\Gamma_f\subset U\times U'$ (the graph of $f$) such that the restriction
$\pr_1|_{\Gamma_f}:\Gamma_f\to U$ of the natural projection $\pr_1:U\times U'\to U$ is 
generically one to one. Generically here means outside of a proper analytic set $I_f$, 
which is called the indeterminacy set of $f$. From the condition of the 
irreducibility of $\Gamma_f$ and Remmert proper image theorem, it follows that $\codim I_f
\ge 2$; in the case of $n=2$ this means a discrete set of points. When $U'=\pp^N$ then 
such $f$ can be defined by an $N$-tuple of meromorphic
functions $f_1,...,f_N$ and $I_f$ is the union of the indeterminacy sets of all $f_k$.
Outside of divisors of poles of $f_k$ we have a mapping to $\cc^N$ viewed as the 
standard affine part of $\pp^N$. With some abuse of notations we shall often write 
$f:U\to \cc^N$ instead of $f:U\to \pp^N$.

\smallskip Our standing assumption about the meromorphic map $f:U\to U'$ is that 
$I_f\subset \{0\}$, \ie either $I_f=\{0\}$ or is empty, and that
\begin{equation}
\eqqno(m-m-prim)
f|_M[M]\subset K,
\end{equation}
for some compact $K$ in $U'$. This will be written as $f:(U,M) \to 
(U', K)$. Let us make precise that $f|_M[M]$ is defined as
\begin{equation}
\eqqno(f-m-m)
f|_M[M] \deff \overline{f(M\setminus \{0\})}.
\end{equation}

This notation is coherent with that from algebraic geometry. In the latter case
$M$ is a complex subvariety of $U$ and $f|_M[M]$ is the so called proper image of $M$
under $f$. 

\newprg[INTER.fam]{Families of analytic disks}

Recall that by an analytic disk in a complex manifold/space $U$ one calls a holomorphic map 
$\phi : \Delta \to U$ of the unit disk $\Delta$ to $U$. The image $\phi (\Delta)$ of an
analytic disk we shall denote as $C$. 
For the purpose of this paper we shall need a precise information about 
intersections of certain holomorphic $1$-parameter families of analytic disks 
$\{C_{\lambda}\}$, all passing through the origin, with our hypersurface $M$. 
Here $\lambda$ is a complex parameter varying in a neighborhood of some 
$\lambda_0\in \cc$. First we shall describe what families $\{C_{\lambda}\}$ 
will occur below.

\smallskip Let $\pi : \hat U\to U$ be a tree of blowings-up over $0\in U\subset \cc^2$.
By saying this we mean that $\pi$ is a composition of a finite number of 
$\sigma $-processes, \ie  $\pi = \pi_1\circ ...\circ \pi_N$. In more details denote by $E_0$ the 
line $\{z_2=0\}$ in the initial neighborhood of zero $U^{(0)}\deff U\subset \cc^2$ with 
``natural coordinates'' $z_1^{(0)} = z_1, z_2^{(0)} = z_2$. In the natural affine charts 
$(U_1^{(1)},z_1^{(1)},z_2^{(1)})$ and $(U_2^{(1)},z_1^{(1)},z_2^{(1)})$ on the first 
blow-up $U^{(1)}$ the blow-down map $\pi_1: U^{(1)}\to U^{(0)}$ can be written as 
\begin{equation}
\begin{cases}
z_1^{(0)} = z_1^{(1)}\cr
z_2^{(0)} = z_1^{(1)}z_2^{(1)}
\end{cases}
\text{ in } \quad U_1^{(1)} \quad \text{ and } \quad
\begin{cases}
z_1^{(0)} = z_1^{(1)}z_2^{(1)}\cr
z_2^{(0)} = z_2^{(1)}
\end{cases}
\text{ in } \quad U_2^{(1)}.
\end{equation}
We denote coordinates both in $U_1^{(1)}$ and $U_2^{(1)}$ with the same letters and call them the
``natural coordinates'' in $U^{(1)}$, this will not lead us to a confusion. The exceptional curve $E_1$ 
is given by $\{z_1^{(1)}=0\}$ in $U_1^{(1)}$ and by $\{z_2^{(1)}=0\}$ in $U_2^{(1)}$. After that one 
blows-up some point on $E_1$ and denotes this as $\pi_2:U^{(2)}\to U^{(1)}$, and so on. Each time in the 
similar way one obtains ``natural coordinates'' 
$z_1^{(k)}, z_2^{(k)}$ on affine charts $U_1^{(k)}$ and $U_2^{(k)}$ of $U^{(k)}$. The exceptional curve 
$E_k$ of $\pi_2:U^{(k)}\to U^{(k-1)}$ writes in these coordinates as either $\{z_1^{(k)}=0\}$ or  as 
$\{z_2^{(k)}=0\}$. We shall not distinguish between $E_k$ and its strict transforms under the further 
$\sigma$-processes, \ie $\pi_{k+1}^*E_k$ will be denoted still as $E_k$ and so on. Therefore the 
exceptional divisor $E$ of the resulting $\pi : U^{(N)}\to U$ represents as the union $E=E_1\cup ...
\cup E_N$. 

\smallskip  For a fixed $k$ between $1$ and $N$ consider the following families:

\begin{equation}
\eqqno(hol-fam1)
\Delta_{\lambda} = \{z_2^{(k)}=\lambda , |z_1^{(k)}| < 1\}_{|\lambda - \lambda_0|<\eps} \text{ if 
the equation of } E_k \text{ is } \{z_1^{(k)}=0\},
\end{equation}
or
\begin{equation}
\eqqno(hol-fam2)
\Delta_{\lambda} = \{z_1^{(k)}=\lambda , |z_2^{(k)}| <1\}_{|\lambda - \lambda_0|<\eps} \text{ if 
the equation of } E_k \text{ is } \{z_2^{(k)}=0\}.
\end{equation}

\medskip\noindent Here $\lambda_0\in \cc$ and $\eps>0$ are chosen arbitrarily. In other words these are 
the families of  disks which intersect $E_k$ orthogonally in 
``natural coordinates'' on $U^{(k)}$. The holomorphic $1$-parameter families we shall be interested in are 
the push-downs of $\Delta_{\lambda}$-s under the blow-down map $\pi_1\circ ... \circ \pi_k:U^{(k)}\to U^{(0)}$.
And they will be denoted actually as $C_{\lambda}$, \ie $C_{\lambda} \deff \pi_1\circ ... \circ 
\pi_k(\Delta_{\lambda})$. The number $1\le k\le N$ will be clear from the context.

\smallskip Families $\{C_{\lambda}\}$ are quite special. Their equations are polynomial in $z_1,z_2$ and 
$\lambda$, this can be easily proved by the induction on the number $N$ of $\sigma$-processes in $\pi$.
In our applications we shall see below that without loss of generality we can assume that the center $\lambda_0$ 
can be taken generic (e.g., avoiding points of intersection of $E_j$ with $E_k$) and $\eps >0$ as small as we wish. 
Generically here means outside of a finite set. Therefore we can assume that 

\smallskip 
\begin{itemize}

\smallskip\item $C_{\lambda}$ \text{ do not degenerate to a point for any } $\lambda$ ;
 
\smallskip\item \text{the tangent cone to any of } $C_{\lambda}$ at zero does not  contain 
either the line \{$z_1=0$\} or  the line \{$z_1=0$\}. 
\end{itemize}

The first assertion is obvious since the proper transform of $C_{\lambda}$ under $\pi$ is 
$\Delta_{\lambda}$, and the latter is a non-constant analytic disk. The second is a bit more subtle.
Would the tangent cone to $C_{\lambda}$ contain more than one line then $\pi_1^*C_{\lambda}$ 
would intersect $E_1$ by more than one point. But this contradicts to the construction of the
family $\{C_{\lambda}\}$. 

\smallskip In what follows we suppose that the line $\{z_1=0\}$ is not in the tangent cone of $C_{\lambda}$ at zero,
the case with $\{z_2=0\}$ can be treated analogously. By genericity we can assume that this holds for 
all $|\lambda -\lambda_0| <\eps$. Therefore there exists a bidisk $\Delta^2(\delta) = 
\Delta (\delta_1)\times \Delta (\delta_2)$ independent of $\lambda$ such that $C_{\lambda}\cap\Delta^2(\delta)$ has 
as its defining function the Weierstrass polynomial 
\begin{equation}
\eqqno(weier)
W_{\lambda}(z) = z_2^d + a^1_{\lambda}(z_1) z_2^{d-1} + ... + a^d_{\lambda}(z_1) \quad\text{ for some } 
\quad d\ge 1,
\end{equation}
where coefficients $a^j_{\lambda}(z_1)$ are holomorphic  both in $z_1$ and in $\lambda$, 
and $a^j_{\lambda}(0)\equiv 0$ again by genericity of $\lambda_0$. Moreover the degree $d$ is independent of 
$\lambda$.  Since 
$C_{\lambda}$ are obviously irreducible at the origin the polynomials $W_{\lambda}$ should be irreducible too. By 
$D_{\lambda }(z_1)$ we denote the discriminant of $W_{\lambda}$. $D_{\lambda }(z_1)$ is holomorphic in 
both variables and $D_{\lambda }(0)\equiv 0$ if $d>1$ because in this case $W_{\lambda}(0,z_2)=z_2^d$
has zero as a root of higher order. Write 
\begin{equation}
\eqqno(discr1)
D_{\lambda }(z_1) = b_k(\lambda )z_1^k + b_{k+1}(\lambda)z_1^{k+1} + ...,
\end{equation}
where $k\ge 1$ and $b_k(\lambda) \not\equiv 0$. Perturbing $\lambda_0$, \ie taking it generically, and taking 
$\eps$ smaller we can assume that the discriminants $D_{\lambda }(z_1)$ of the equations $\{W_{\lambda } =0\}$ 
of $C_{\lambda}$ do not vanish for $z_1\in \Delta (\delta_1)\setminus 0$ for some $\delta_1 >0$ independent of 
$\lambda$. Therefore in the bidisk $\Delta^2 (\delta)$ our curves $C_{\lambda}$ are given by the 
equations
\begin{equation}
\eqqno(curves1)
C_{\lambda } = \left\{z_2 = h_{\lambda}(z_1^{\frac{1}{d}})\right\}, h_{\lambda}(0)=0,
\end{equation}
see \cite{F}. Dependence of $h_{\lambda}$ on $\lambda$ stays to be holomorphic. If $d=1$ one readily gets the same form 
\eqqref(curves1) for $C_{\lambda}$.

\newprg[INTER.int]{Intersection of real hypersurfaces with families of analytic disks}

Let $R$ be the intersection of a non-constant analytic disk $C = \phi (\Delta)\ni 0$ with the 
strictly pseudoconvex real analytic hypersurface $0\in M\subset\cc^2$ as above, \ie $R = 
\phi^{-1}(C\cap M)$. Then $R$ is a one dimensional real analytic set, and it has a non-empty 
(one dimensional) interior provided $C\cap U^+\not= \emptyset$. More precisely  $R = S\cup \Gamma$, 
where $S$ is a discrete in $\Delta$ set of points $\{s_k\}$ and $\Gamma$ is a locally finite union  
of smooth arcs $\{\gamma_l\}$ with ends on $\{s_k\}$.

\smallskip Now consider a holomorphic $1$-parameter family $\{\phi_{\lambda} : |\lambda - \lambda_0|
<\eps \}$ of analytic disks $\phi_{\lambda} : (\Delta , 0)\to (\cc^2, 0)$ such that $C_{\lambda}
=\phi_{\lambda}(\Delta)\cap \Delta^2(\delta) = \{ (z_1, h_{\lambda}(z_1^{\frac{1}{d}}):
z_1\in \Delta_{\delta_1}\}$ are as above with $R_{\lambda}$, $\Gamma_{\lambda}$ and $S_{\lambda}$
having an obvious meaning in the parameter case.

\begin {lem}
\label{intersect1}
For $\eps >0$ small enough the following holds:

\smallskip\sli either for all $|\lambda - \lambda_0| <\eps $ the $1$-dimensional part 
$\Gamma_{\lambda}$ of $R_{\lambda}:=\phi_{\lambda}^{-1}(M\cap C_{\lambda})$ contains a 

\quad component $\gamma_{\lambda}$ with an end point at zero;

\smallskip\slii or  $R_{\lambda} =\{0\}$ for all $\lambda$ and then $C_{\lambda}\setminus 0 \subset U^-$;

\smallskip\sliii or  there exist
sequences $\lambda_n\to \lambda_0$ and $\eps_n\to 0$ 
such that (\sli holds for all $|\lambda -\lambda_n|<\eps_n$ 

\quad and for all $n$.
\end {lem}
\proof Equations of analytic disks $C_{\lambda}$ and hypersurface $M$ can be written as  

\smallskip
a) $C_{\lambda} = \{z_2 = h_{\lambda}(z_1^{\frac{1}{d}})\}$, $h_{\lambda}(0)= 0$ as in \eqqref(curves1).
 
\smallskip b) $M = \{y_2 = |z_1|^2 + \sum_{k, l\ge 2} a_{k\bar l}(x_2)z_1^k\bar z_1^l\}$ as in \eqqref(eqm-2).
 
\smallskip\noindent Making the substitution 

\begin{equation}
\begin{cases}
z_1\to z_1^d \cr
z_2\to z_2
\end{cases}
\end{equation}
we rewrite (a) and (b) in the form 

\begin{equation}
\eqqno(phi-and-m)
\begin{cases}
C_{\lambda} = \{z_2 = z_1^qh_{\lambda}(z_1) \} \text{ with some } q\ge 1, \cr
M = \{y_2 = |z_1|^{2d} + \sum_{k, l\ge 2} a_{k\bar l}(x_2)z_1^{dk}\bar z_1^{dl}\},
\end{cases}
\end{equation}
with some (other) holomorphic $h_{\lambda}$ such that $h_{\lambda }(0)\not\equiv0$ as a function of $\lambda$.

\medskip\noindent{\slsf Case 1. } {$q>2d$}. In that case we obviously and directly 
from \eqqref(phi-and-m) get that $C_{\lambda}\cap \bar U^+ = \{0\}$ for all $\lambda$, 
\ie we get the option (\slii of the lemma.

\medskip\noindent{\slsf Case 2.} {$h_{\lambda_0}(0)\not=0$.} Taking $\eps >0$ smaller, if necessary, we can 
suppose that $h_{\lambda }(0) \not=0$ for all $|\lambda -\lambda_0|<\eps$. Write 
$h_{\lambda}(0)=a(\lambda)e^{i\theta (\lambda)}$.  Using polar coordinates $z_1 = r_1e^{\i\phi_1}$ 
and $z_2 = r_2e^{i\phi_2}$ we get that on $C_{\lambda}\cap M$ one has 

\begin{equation}
\eqqno(cases)
\begin{cases}
r_2\sin \phi_2 = r_1^{2d} + O(r_1^{4d}), \cr
r_2\sin \phi_2 = a(\lambda)r_1^q\sin (q\phi_1 + \theta (\lambda)) + O(r_1^{q+1}), \cr
r_2\cos \phi_2 = a(\lambda)r_1^q\cos (q\phi_1 + \theta (\lambda)) + O(r_1^{q+1}).
\end{cases}
\end{equation}
Now we have the following two subcases.

\smallskip\noindent{\slsf Subcase 2a.} {$q<2d$.} From \eqqref(cases) we get that 
\[
r_1^{2d} + O(r_1^{4d}) = a(\lambda)r_1^{q} \sin (q\phi_1 + \theta (\lambda)) + O(r_1^{q+1}),
\]
which implies 
\[
a(\lambda)\sin (q\phi_1 + \theta (\lambda)) = r_1^{2d-q} + O(r_1^{4d-q}) + O(r_1) = O(r_1).
\]
The last equation has $2q$ curves of solutions as $r_1$ tends to zero:
\[
\phi_1 = -\frac{1}{q}\theta (\lambda) + \frac{1}{q}\arcsin \left(\frac{1}{a(\lambda)}O(r_1)\right)
+ \frac{2\pi k}{q}, \quad -q\le k\le q-1,
\]
and all these curves end at the origin. 

\smallskip After we had determined $\phi_1$ as a function of $r_1$
from equations one and two of \eqqref(cases) we can easily find $z_2 = r_2(\cos \phi_2 + i\sin\phi_2)$
as a function of $\r_1$ from the equations one and three of \eqqref(cases). Since $q\ge 1$
we shall have that $z_2(r_1)\to 0$ as $r_1\to 0$. I.e., the option (\sli of our lemma realizes.

\smallskip\noindent{\slsf Subcase 2b.} {$q=2d$.} Again, from \eqqref(cases) we get that 
\[
r_1^{2d} + O(r_1^{4d}) = a(\lambda) r_1^{2d} \sin (q\phi_1 + \theta (\lambda)) + O(r_1^{2d+1}),
\]
which implies 
\begin{equation}
\eqqno(q=2p)
a(\lambda)\sin (q\phi_1 + \theta (\lambda)) = 1 + O(r_1).
\end{equation}
If $a(\lambda_0)<1$ then equation \eqqref(q=2p) has no solutions for small $r_1$ and we fall to the
option (\slii of our lemma. If $a(\lambda_0)>1$ then we are again in  (\sli with the curves of solutions

\[
\phi_1 = -\frac{1}{q}\theta (\lambda) + \frac{1}{q}\arcsin \left(\frac{1}{a(\lambda)}(1 + O(r_1))\right)
+ \frac{2\pi k}{q}, \quad -q\le k\le q-1.
\]
Finally, if $a=1$ 
then by open mapping theorem one finds $\lambda $ arbitrarily close to $\lambda _0$ such that 
$a(\lambda) <1$ and the one gets the option (\sliii of our lemma.

\smallskip\noindent{\slsf Case 3. }  If $h_{\lambda_0 }(0)=0$
we can take $\lambda_n$ arbitrarily close to $\lambda_0$ such that $h_{\lambda_n }(0)\not=0$ but 
$a(\lambda_n) \deff |h(\lambda_n)|$ small. We fall into the assumptions of Case 2 with an additional 
condition that $a(\lambda )<1$. This gives us the option (\sliii of our lemma.

\smallskip\qed

\newsect[STR]{Strict transform of real analytic hypersurfaces under a modification}

\newprg[STR.str]{Strict transform}

Let $M\ni 0$ be a real analytic strictly pseudoconvex hypersurface near the origin
in $\cc^2$. Let $\pi : \hat U\to U$ be a tree of blowings-up over the origin, $U$
stands here for some neighborhood of zero in $\cc^2$. The exceptional divisor of 
$\pi$ is denoted as $E$. Denote by 
\begin{equation}
\eqqno(str-tr1)
\pi^*M \deff \overline{\pi^{-1}(M\setminus 0)}
\end{equation}
the proper preimage (or  strict transform) of $M$ under $\pi$. Set $M_0^*\deff \pi^*M\cap E$. 
Note that $M^*_0$ is connected. This follows from the fact that this set is the intersection 
of connected sets 
\begin{equation}
\eqqno(str-tr2) 
M_0^* = \bigcap_{\eps > 0} \overline{\pi^{-1}\left((B^4_{\eps}\setminus 0)\cap M\right)},
\end{equation}
here $B^4_{\eps}$ stands for the ball of radius $\eps$ centered at the origin of $\cc^2$. 
We are going to prove that $M_0^*$ is a {\slsf massive set} in the topology of $E$.

\begin{lem}
\label{massive-m}
Let $\hat p\in M_0^*$ be a point on the strict transform of $M$ over the origin which is 
not a singular point of $E$ and let $E_k$ be an irreducible component of $E$ containing 
$\hat p$. Then $\hat p$ is an accumulation point of $int (M_0^*\cap E_k)$. In particular 
$M^*_0\cap E_k$ has the non-empty interior in the topology of $E_k$ for every $E_k$ such 
that $M^*_0\cap E_k \not\subset \Sing E$.
\end{lem}
\proof Let $\hat p\in M^*_0\cap E_k$ be as in this lemma. From \eqqref(str-tr2) it is 
clear that we can find a sequence $\hat p_n\to \hat p$
such that $p_n\deff \pi (\hat p_n)\in U^+$ and then necessarily $p_n\to 0$. Moreover, since
$U^+$ is open we can choose these $\hat p_n$ generically. For any fixed $n>> 1$ find 
natural coordinates $z_1^{(k)},z_2^{(k)}$ in an affine neighborhood of  $\hat p_n$ and set 
$\lambda_0 \deff z_1^{(k)}(\hat p_n)$ or  $\lambda_0 \deff z_2^{(k)}(\hat p_n)$ depending
on what is the equation of $E_k$ in these coordinates, see \eqqref(hol-fam1) and 
\eqqref(hol-fam2). Let $\{\Delta_{\lambda}\}_{|\lambda -\lambda_0|<\eps}$ and 
$\{C_{\lambda}\}_{|\lambda -\lambda_0|<\eps}$ be the holomorphic $1$-parameter families
constructed as there. By the genericity of the choice of $\hat p_n$ these families
can be chosen generically as well. Since these disks cut $U^+$ the case (\slii of Lemma 
\ref{intersect1} does not occur and we conclude that either $\lambda_0$ is an interior 
point of $M_0^*$ or it can be approximated by interior points. This gives in its turn 
the approximation of $\hat p$ by the interior points of $M_0^*$.

\smallskip Indeed, suppose we are under the case (\sli of that lemma. Then $C_{\lambda}$ intersects
$M$ by an $1$-dimensional real analytic set  which accumulates to zero. Let $\gamma_{\lambda}$
be some $1$-dimensional local component of this set which accumulates to $0$. Its strict
transform $\hat \gamma_{\lambda}\deff \overline{\pi^{-1}(\gamma_{\lambda}\setminus 0)}$ accumulates 
$\lambda \in E_k$, \ie $\lambda$ is viewed as the point of intersection of $\Delta_{\lambda}$
with $E_k$. Therefore $\lambda$ belongs to $M^*_0\cap E_k$ for all $\lambda $ close to 
$\lambda_0$. The case (\sliii is clear.

\smallskip\qed

\begin{rema} \rm
\label{fat-m}
Observe that $E_k\cap \Sing E$ consists from a finite set of points. From connectivity of $M_0^*\cap E_k$
it follows that if $M_0^*\cap E_k$ is non-empty and is contained in $\Sing E$ then $M_0^*\cap E_k = 
\{\text{point}\}$.
\end{rema}

\newprg[CONT.transf]{Strict transform to the first two blowings-up: example.}

It is not necessary for the proofs of this paper but is very instructive to 
compute the strict transform of a real analytic hypersurface onto first
few blowings-up. 

\smallskip\noindent 1) Write the equation \eqqref(eqm-2) in the form 
\begin{equation} 
\eqqno(eqm-3)
M = \{\im z_2 =  |z_1|^2 + r(z,\bar z)\},
\end{equation}
If we denote by $M_1$ the proper transform of $M$ under $\pi_1$ and use notations of 
subsection \ref{INTER.fam} we see that $M_1\setminus E_1$ has equations

\begin{equation}
\eqqno(step11)
M_1\setminus E_1  = \{\im (z_1^{(1)}z_2^{(1)}) = |z_1^{(1)}|^2 + r_1^{(1)}(z^{(1)},\bar z^{(1)})\}  
\quad\text{ in }\quad U_1^{(1)},
\end{equation}

\noindent where $r_1^{(1)}(z^{(1)},\bar z^{(1)}) = r(z_1^{(1)},z_1^{(1)}z_2^{(1)}, 
\bar z_1^{(1)}, \bar z_1^{(1)}\bar z_2^{(1)})$, and 

\begin{equation}
\eqqno(step12)
M_1\setminus E_1 = \big\{\im z_2^{(1)} = |z_1^{(1)}|^2|z_2^{(1)}|^2 + r_2^{(1)}(z^{(1)},\bar z^{(1)})\big\} 
\quad\text{ in }\quad U_2^{(1)},
\end{equation}
where $r_2^{(1)}(z^{(1)},\bar z^{(1)}) = r(z_1^{(1)}z_2^{(1)}, z_2^{(1)}, \bar z_1^{(1)}\bar z_2^{(1)}, 
\bar z_2^{(1)})$. We see that the rests $r_1^{(1)}$ and $r_2^{(1)}$ have order of vanishing at zero not 
less than the original $r$.

\smallskip The closure  of $M_1\setminus E_1$ in $U_2^{(1)}$, that is actually 
$M_1\cap U_2^{(1)}$, is a smooth hypersurface (one easily checks that the gradient never 
vanishes) with the same equation as \eqqref(step12) and $M_1\cap U_2^{(1)}\supset 
(E_1\cap U_2^{(1)})$. Now remark that $M_1\setminus E_1$ in $U_1^{(1)}\setminus E_1$ is defined 
by
\begin{equation}
M_1\setminus E_1 = \left\{\im \left(\frac{z_2^{(1)}}{\bar z_1^{(1)}}\right) = 1 + 
\frac{r_1^{(1)}(z^{(1)},\bar z^{(1)})}{z_1^{(1)}\bar z_1^{(1)}}\right\},
\end{equation}
where $\frac{r_1^{(1)}(z^{(1)},\bar z^{(1)})}{z_1^{(1)}\bar z_1^{(1)}} = O(z^{(1)})$.

\begin{figure}[h]
\centering
\includegraphics[width=1.7in]{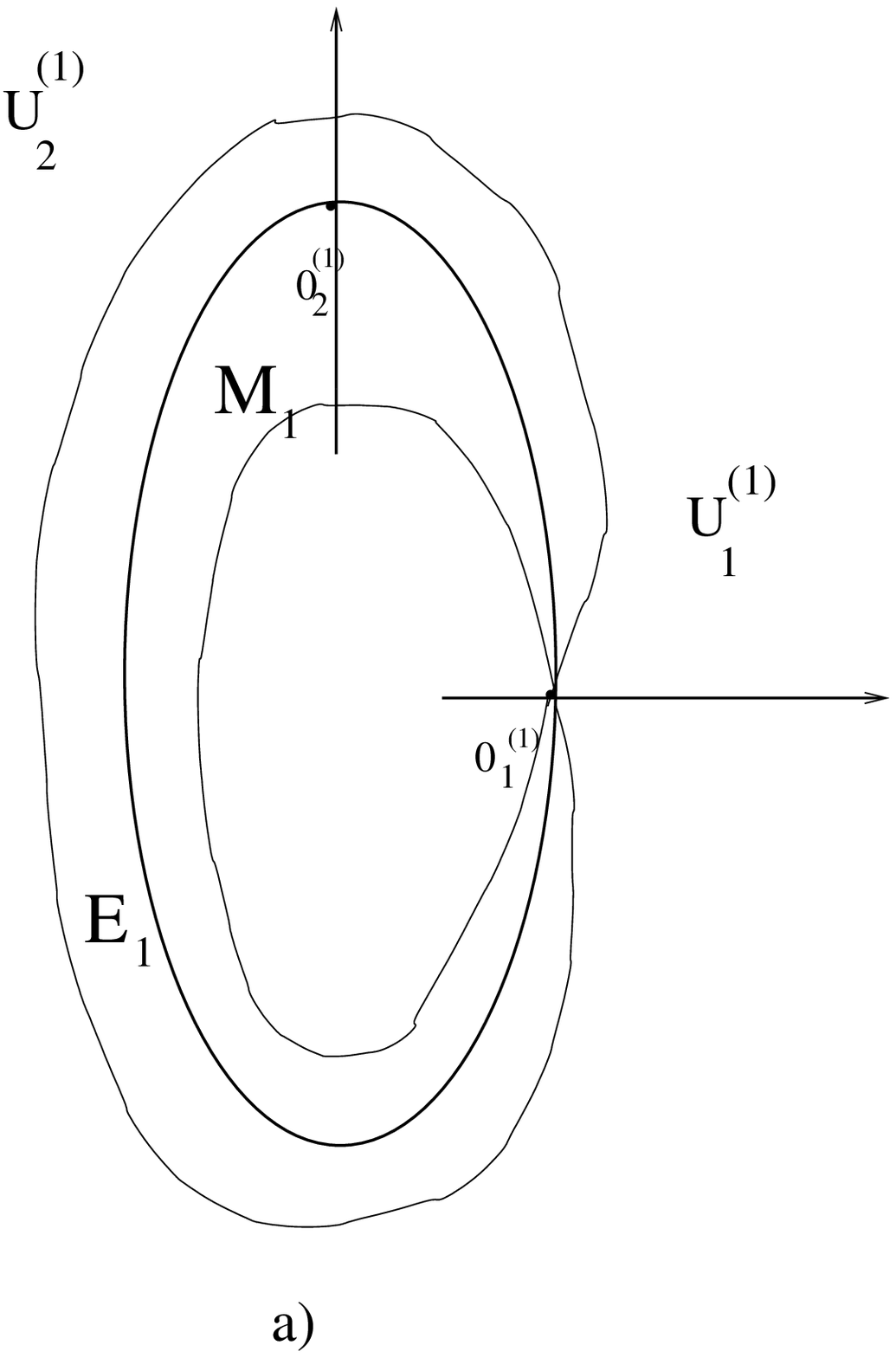}
\quad\quad\quad
\includegraphics[width=1.6in]{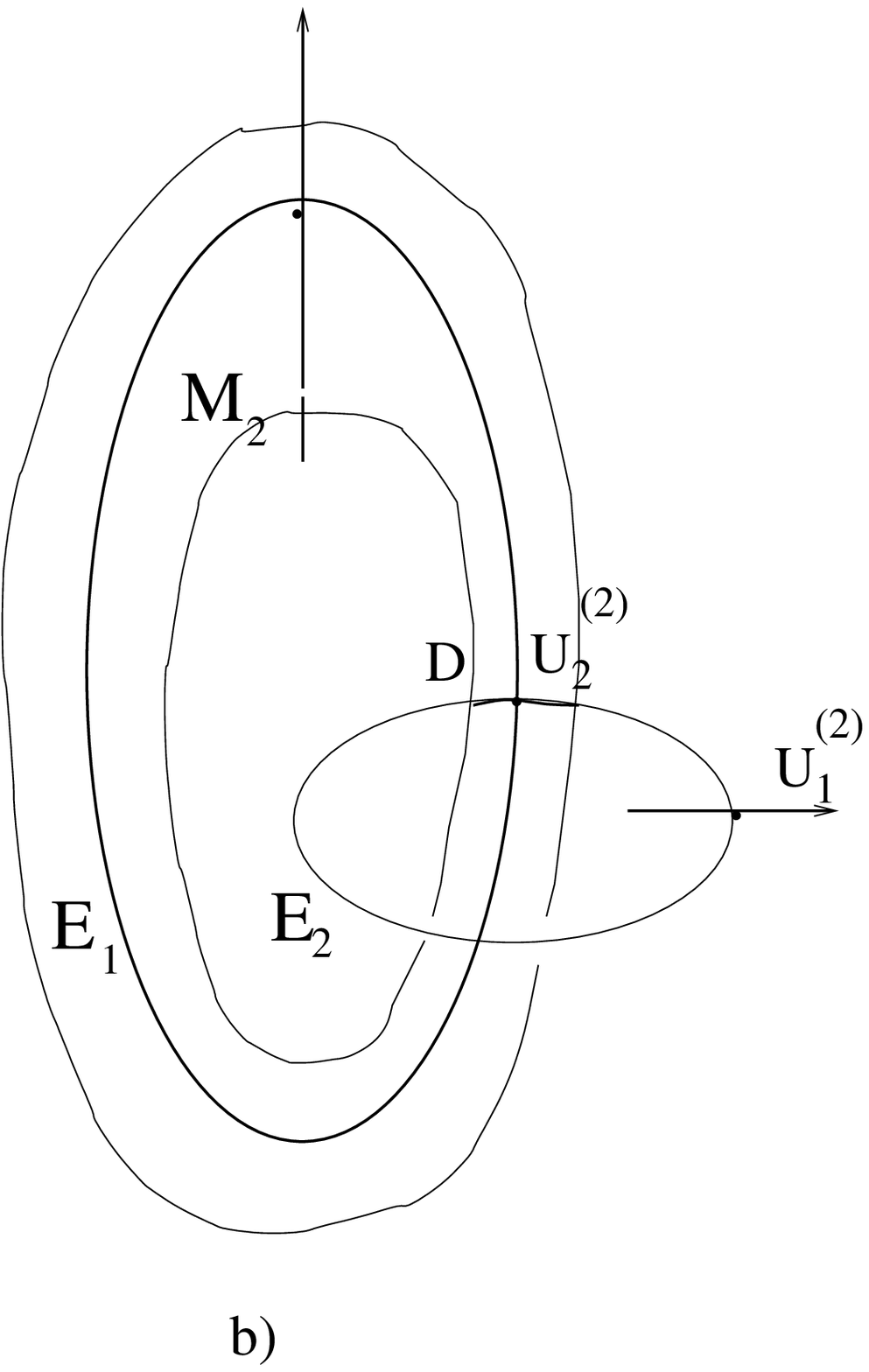}
\caption{The first proper transform $M_1$ of $M$, see picture (a), contains the
first exceptional curve $E_1$ entirely. The second $M_2$, see picture (b), contains
$E_1$ and only a part of $E_2$.}
\label{step-1}
\end{figure}

\smallskip We see that $M_1\cap U_1^{(1)}$ is a real cone with vertex at the origin. On 
the diagram of moduli $(r_1,r_2) = (|z_1^{(1)}|,|z_2^{(1)}|)$ it is tangent to the cone $r_2
\ge r_1$, see Fig. \ref{step-1} (a). At this stage the proper 
transform $M_1$ of $M$ contains the entire exceptional curve $E_1$.

\smallskip\noindent 2) Now let us blow-up the point $0_1^{(1)}$, by which we denote the origin in $U_1^{(1)}$. 
In what follows coordinates $(z_1^{(1)},z_2^{(1)})$ will be redefined simply as $(z_1,z_2)$ in order 
to simplify the notations, the same for $z^{(2)}$ below. We get the following equations in charts 
$U_1^{(2)}$ and $U_2^{(2)}$:

\begin{equation}
\eqqno(step21)
\im (z_1^2z_2) = |z_1|^2  + r_1(z, \bar z)\quad\text{ in }\quad U_1^{(2)},
\end{equation}
and
\begin{equation}
\eqqno(step22)
\im (z_1z_2^2) = |z_1|^2|z_2|^2 + r_2(z,\bar z) \quad\text{ in } U_2^{(2)},
\end{equation}
with an appropriate $r_1$ and $r_2$. 
If we denote by $M_2$ the proper transform of $M_1$ under $\pi_2:U^{(2)}\to U^{(1)}$ and by $E_2$
the corresponding exceptional curve then for $M_2\setminus E_2$ in $U_1^{(2)}\setminus E_2$ 
and for $M_2\setminus E_2$ in $U_2^{(2)}\setminus E_2$ we get correspondingly the 
equations
\begin{equation}
\im \left(\frac{z_1}{\bar z_1}z_2\right) = 1 + O(z)\quad \text{ and } \quad
\im \left(\frac{z_2}{\bar z_2}\cdot \frac{1}{\bar z_1}\right) = 1 + O(z).
\end{equation}
Now we see that $M_2$ intersects the second exceptional curve $E_2$  by the 
closed disk $D$, which in coordinates of $U_1^{(2)}$, corr. of $U_2^{(2)}$, is given as
\begin{equation}
\eqqno(eq-of-d)
D = \{z_1= 0, |z_2|\ge 1\} \quad\text{ corr. as }\quad D = \{z_2 = 0, |z_1|\le 1\},
\end{equation}
see Fig. \ref{step-1} b). 

\newsect[CONT]{Continuity of the mapping}

\newprg[CONT.lim]{The limit set of $f$ from the concave side of $M$}

Denote by 

\begin{equation}
\eqqno(lim1)
\lim_{\substack{z\to 0\\z\in U^+}} f(z)  \quad\text{ correspondingly } \quad 
\lim_{\substack{z\to 0\\z\in M\setminus 0}} f(z),
\end{equation}
the sets of cluster points of all sequences $f(z_n)$ when $z_n\in U^+$, $z_n\to 0$ (corr. when
$z_n\in M\setminus 0$, $z_n\to 0$).

\begin{lem}
\label{limit-f}
Let a meromorphic mapping $f$ satisfies \eqqref(m-m-prim), where $M$ is as above and  $M'\comp
U'$ is an arbitrary compact. Then

\begin{equation}
\eqqno(lim2)
\lim_{\substack{z\to 0\\ z\in U^+}} f(z) = 
\lim_{\substack{z\to 0\\ z\in M\setminus 0}} f(z) .
\end{equation}
\end{lem}
\proof Let $\pi:\hat U\to U$ be a resolution of indeterminacies of $f$, \ie a tree of blowings-up 
over the origin, such that the pull-back $\hat f \deff f\circ \pi$ of $f$ to $\hat U$ is holomorphic.
Furthermore, let $z_n\in U^+$ be a sequence of points such that $z_n\to 0$ and $f(z_n) \to w_{\infty}$. 
The latter is some point of $M'$.
Denote by $\hat z_n\in \hat U^+\deff \pi^{-1}(U^+)$ the preimages of $z_n$ under $\pi$. Taking a 
subsequence we can suppose that $\hat z_n$ converge to some $\hat z_{\infty}\in E\deff \pi^{-1}(0)$, 
the latter is the exceptional divisor of the modification $\pi$. We have that $\hat f(\hat z_n) \to  
w_{\infty}= \hat f(\hat z_{\infty})$. Take the family $\{\Delta_{\lambda}\}_{|\lambda -\lambda_0|}<\eps$
of disks as in \eqqref(hol-fam1) and \eqqref(hol-fam2) in such a way that $\Delta_{\lambda_0}\cap E =
\hat z_{\infty}$. 

\smallskip Before applying Lemma \ref{intersect1} to the corresponding family $C_{\lambda}$ we perturb 
$\lambda_0$ and take $\eps >0$ in order to fulfill the usual assumptions we impose on our family.
Remark also that we can find $\lambda$-s arbitrarily close to $\lambda_0$ in such a way that in addition 
$C_{\lambda}\setminus 0$ cannot be contained in $U^-$. For example take $\lambda$ such that $\Delta_{\lambda}$
contains some $\hat z_n$, this means that $C_{ \lambda} \ni z_n\in U^+$ (or  any $\lambda$ close to this). 
Therefore the case (\slii of Lemma \ref{intersect1} will not occur. Remark that we perturbed also ours
$\hat z_{\infty}$ and $w_{\infty}$.

\smallskip\noindent{\slsf Case 1}. If we are in the case (\sli of Lemma \ref{intersect1} then 
\begin{equation}
\eqqno(lim3)
\lim_{\substack{z\to 0\\ z\in C_{\lambda_0}}} f(z) = \lim_{\substack{\hat z\to \hat z_{\infty}\\ 
\hat z\in \Delta_{\lambda_0}}} \hat f(z) = w_{\infty}.
\end{equation}
But $C_{\lambda_0}\cap M$  contains a real analytic curve $\gamma$ accumulating 
to zero in this case. The limit in \eqqref(lim3) is the same as the limit along this $\gamma$
which is a subset of $M$. This proves the inclusion 
\begin{equation}
\eqqno(lim4)
\lim_{\substack{z\to 0\\ z\in U^+}} f(z) \subset
\lim_{\substack{z\to 0\\ z\in M\setminus 0}} f(z)
\end{equation}
in this case. 

\smallskip\noindent{\slsf Case 2.} Suppose we are in the case (\sliii of Lemma \ref{intersect1}. Let 
$\lambda_n\to \lambda_0$ be
as there and let $\hat z_n'$ be the point of intersection of $\Delta_{\lambda_n}$ with $E$. Remark that
$\hat z_n'\to \hat z_{\infty}$. Let $\gamma_n$ be a component of $C_{\lambda_n}\cap M$ accumulating to
zero. Then 
\begin{equation}
\eqqno(lim7)
\hat f(\hat z_n') = \lim_{\substack{\hat z\to \hat z_n'\\ \hat z\in \Delta_{\lambda_n}}} \hat f(\hat z) =
\lim_{\substack{z\to 0\\ z\in C_{\lambda_n}}} f(z) = \lim_{\substack{z\to 0\\ z\in \gamma_n}} f(z)
\subset \lim_{\substack{z\to 0\\ z\in M}} f(z).
\end{equation}
By holomorphicity of $\hat f$ we have that $\hat f(\hat z_n') \to \hat f(\hat z_{\infty}) = w_{\infty}$. 
I.e., the inclusion \eqqref(lim4) is proved also in this case.

\smallskip In fact up to now we proved that for some $\hat z_n^{''}\in E, \hat z_n^{''}\to \hat z_{\infty}$
we have that $\hat f(\hat z_n^{''})\in f|_M[M]$. But the considerations applied in the Case 2 show that
this implies that $\hat f(\hat z_{\infty})\in f|_M[M]$ as well. I.e., the inclusion \eqqref(lim4) 
is proved.

\smallskip The inverse inclusion is obvious. For if 
\[
w_{\infty }= \lim_{\substack{z_n\to 0\\z_n\in M\setminus0}} f(z_n),
\]
then we find $\tilde z_n\in U^+$ sufficiently close to $z_n$ and the $\lim f(\tilde z_n)$ will be the same.

\smallskip\qed

\begin{rema} \rm 
The set in the right hand side of \eqqref(lim2) we shall denote as $f|_M[0]$. Observe
the obvious inclusion 
\begin{equation}
\eqqno(lim5)
\lim_{\substack{z\to 0\\ z\in M\setminus 0}} f(z) = f|_M[0] \subset f|_M[M].
\end{equation}
\end{rema}

\newprg[CONT.proof]{Continuity of the mapping from the concave side of $M$}

Now we can proof the announced continuity of $f$ from the side $U^+$ of $M$.

\begin{lem}
\label{cont-f}
In the conditions of Lemma \ref{limit-f} suppose, in addition, that the compact $M'$ 
does not contain any germ of a complex curve. Then the restriction $f|_{\bar U^+\setminus 0}$ 
is continuous up to $M$.
\end{lem}
\proof We need to prove the continuity at zero only. Moreover, due to Lemma 
\ref{limit-f} the only thing to prove is that 
\begin{equation}
\eqqno(lim6)
f|_M[0] = \{ point\}.
\end{equation}
Let $\pi : \hat U\to U$ be, as above, the resolution of indeterminacies of $f$ and $E$
its exceptional divisor.   It is clear that 
\begin{equation}
\eqqno(f-at-0)
f|_M[0] = \hat f\left(M_0^*\right).
\end{equation}
By $E^{\const}_{\hat f}$ denote the union of components of $E$ on which the lift $\hat f$ 
is constant. Therefore if $M_0^*\subset E^{\const}_{\hat f}$ then the 
lemma is proved. Suppose that this is not the case. This means that there
exists an irreducible component $E_k$ of $E$ such that $M_0^*\cap E_k\not=
\emptyset $ and $\hat f|_{E_k} \not\equiv\const$. If $M_0^*\cap E_k\subset \Sing E$, 
\ie is just one point $\hat p$, see Remark \ref{fat-m}, then there should be another component
$E_l$ of $E$ intersecting $E_k$ at this point. $\hat p$ cannot be an isolated point 
of $E_l$ too and we can find a point $\hat q\in M_0^*$ close to $\hat p$ which is 
not singular for $E$ and then Lemma \ref{intersect1} applies to $\hat q$ on the
place of $\hat p$. 

\smallskip According to this lemma there exists $\hat q\in E_0$ close to 
$\hat p$ and a neighborhood $\hat V$ of $\hat q$ on $E_k$ (or  on $E_l$) such that 
$\hat V\subset M^*_0$. But this means that $\hat f|_{E_k}(\hat V) \subset M'$ (or  
$\hat f|_{E_l}(\hat V) \subset M'$). Since $M'$ does not contain germs of complex curves
this means that $\hat f$ is constant on $E_k$ (or  on $E_l$). In the former case this is
a contradiction. In the latter $E_k\cap M_0^* = \{\hat p\}\in E^{\const}_f$. Therefore 
$M_0^*\cap E_k\subset E^{\const}_{\hat f}$ for every component $E_k$ of $E$
which intersects $M_0^*$ and the lemma is proved.

\smallskip\qed

\begin{rema} \rm
If no further assumption are imposed on the image set $K$ mapping $f$ (being 
continuous from the concave side $U^+$ of $M$) can be still {\slsf meromorphic}
in general, \ie it can happen that $0$ is really an indeterminacy point of $f$.
Such an example was constructed in \cite{IM}. 
\end{rema}

\newsect[APL]{Applications and generalizations}

\newprg[APL.proof1-2]{Proof of Corollaries 1, 2 and 3.} 
Let a germ $f:(M,x)\to (M',x')$ of a holomorphic mapping from a real analytic hypersurface 
$M\subset \cc^n$ to a compact locally spherical $M'\subset \cc^{n'}$
be given. By the result of \cite{SV} $f$ meromorphically extends along any given 
$CR$-path $\gamma\subset M$ starting at $x$. All we need to prove is that this extension
is holomorphic. Denote by $0\in M$ the point in which we shall prove the holomorphicity
of $f$.

\smallskip Take a vector $\v\in T_0^cM$ such that $L_M(0)[\v]\not=0$. Taking such $\v$ 
generically we can preserve this condition 
and, additionally, taking  a transverse to $M$ vector $\n$ we can achieve that $L\cap I_f$ 
is discrete for a subspace $L$ of $\cc^n$ spanned by $\v$ and $\n$. Then $L\cap I_f$ 
will be discrete as well. All we need is to prove
our theorem for the restriction of $f$ to this subspace. Indeed,  after a coordinate
change we can suppose that $L = \{z_3=...=z_n=0\}$ and that $f$ is meromorphic in 
the unit polydisk $U=\Delta^n$. Shrinking $U$ if necessary and assuming that our theorem 
is proved when $n=2$ we get that  $f|_{L\cap U}$ is holomorphic and therefore the graph 
$\Gamma|_{f|_{L\cap U}}$ is Stein. Take a Stein neighborhood $V$ of $\Gamma|_{f|_{L\cap U}}$ 
in $U\times \pp^N$. Then for every $z^{''} = (0,0,z_3,...,z_n)$ close to zero we have that 
$\hat f(z) \deff (z, f(z))$ is holomorphic in a neighborhood of $\d\Delta^2\times \{z^{''}\}$
with values in $V$. The holomorphicity of $f$ follows now from the Hartogs extension theorem 
for holomorphic functions. Therefore from now on we shall assume that $n=2$, $M$ strictly 
pseudoconvex at $0$ and zero is the only eventual indeterminacy point of $f$.

\smallskip Let $M'$ be our compact locally spherical hypersurface in $\cc^{n'}$. This 
means that for every point $x'\in M'$ there exists a neighborhood $U'\ni x'$ in $\cc^{n'}$
and a biholomorphism $\Phi$ of $U'$ with values in $\cc^{n'}$ such that $\Phi (U'\cap M')
\subset \sss^{2n'-1}$. In \cite{BS} it was proved that the universal cover $\tilde M'$ 
of $M'$ admits a biholomorphic {\slsf development} mapping $F:\tilde M'\to \sss^{2n'-1}$.
Moreover, construction in \cite{BS} obviously gives a complex neighborhood $\tilde V'\supset 
\tilde M'$ such that the covering map $\pi : \tilde M'\to M'$ extends as a locally 
biholomorphic map to $\tilde V'$ and this $f$ is an embedding of $\tilde V'$ to $\cc^{n'}$.

\smallskip Shrinking $(U,M)$, if necessary, we can suppose that $f(\overline{U^+})\subset U'$, 
where this time $U'$ is a neighborhood of $0'=f(0)\in M'$ such that $\pi : \tilde U'\to U'$ 
is a biholomorphism. Here $\tilde U'$ is an appropriate neighborhood of some preimage $\tilde 0'$ 
of $0'$ by $\pi$ in $\tilde V'$.
As a result the composition $H\deff F\circ\pi^{-1}\circ f:\overline{U^+}\to \cc^{n'}$ is well 
defined and maps $M$ to $\sss^{2n'-1}$. Moreover, it extends holomorphically to a neighborhood
of any point $x\in M\setminus \{0\}$. Applying the quoted theorem of Pinchuk we extend
$H$ holomorphically to a neighborhood of $M$, denote this neighborhood by $U$ again. This gives us the 
desired holomorphicity of $f$ and proves Corollaries \ref{hol-ext1} and \ref{hol-ext2}.

\smallskip\noindent{\slsf Proof of Corollary 3.} As it was said in Introduction the distribution 
of complex tangents on a strictly pseudoconvex hypersurface is contact and by Theorem of 
Gromov, see \cite{G} one can draw a CR-path between any two points of $M$. Since, in 
addition $M$ is supposed to be simply connected our $f$ extends meromorphically to a neighborhood
of $M$ and therefore on $\bar D$. If $n=2$ we can directly apply our Theorem \ref{mer-cont}
and continuity of $f$ up to the boundary. This gives us the statement of Corollary \ref{hol-ext3}.

\begin{rema} \rm
Corollary \ref{hol-ext3}  most probably stays to be true for all $n\ge 2$ in the source and not
only for $n=2$. But the proof requires the study of multidimensional blowings-up of real hypersurfaces
and is out of the range of this paper.
\end{rema}

\newprg[APL.corr]{Proof for meromorphic correspondences}

Let $U$ and $U'$ be complex manifolds. Recall that a meromorphic correspondence 
$f:U\to U'$ is an  analytic subset $\Gamma_f\subset U\times U'$
such that the restriction $\pr_1|_{\Gamma_f}:\Gamma_f\to U$ of the natural projection $\pr_1:U\times
U'\to U$ onto $\Gamma_f$ is proper and generically $d$ to $1$. Here $d\ge 1$ is called the 
order of $f$ and $\Gamma_f$ is its graph. Meromorphic correspondence of order $\adyn$ is just
a meromorphic mapping. In general $f$ is called a $d$-valued meromorphic correspondence. 
If $I_f=\emptyset$ the correspondence is called holomorphic. Here $I_f$ stands for the set of 
indeterminacy points of $f$, \ie $x\in I_f$ if $\dim \pr_1|_{\Gamma_f}^{-1}(x)\ge 1$.

\smallskip Let us prove now Corollary \ref{corr-cont}. The symmetric power $\sym^d(U')$ of $U'$ of order
$d$ is a normal complex space and $f$ naturally defines a meromorphic mapping $f^d:U\to \sym^d(U')$.
Condition that $f(m)\subset K$ for $m\in M\setminus \{0\}$ implies (is equivalent to) that 
$f^d(m)\in  \sym^d(K)$. Would $\sym^d(K)$ contain a germ of a non-constant complex curve 
$\phi^d(\lambda) = \sym (\phi_1(\lambda),...,\phi_d(\lambda))$ then $K$ would contain all $\phi_k$.
At least one of them should be non-constant. Contradiction. Therefore $\sym^d(K)$ does not contain
such a germ. Now one can literally repeat the proof of Theorem \ref{mer-cont} for the meromorphic
mapping $f^d:U\to \sym^d(U')$ to get the conclusion of Corollary \ref{corr-cont}.

\newsect[APP]{Appendix: result of Shafikov-Verma}

In this appendix we explain that the paper \cite{SV} contains the following statement. 

\begin{thm}
\label{SVT}
Let $M$ be a smooth real analytic minimal hypersurface in $\cc^n$ and let
$M'$ be a smooth compact real algebraic hypersurface in $\cc^{n'}, 1 < n \le n'$. 
Then every germ $f:(M,x) \to (M',x')$ of a holomorphic map from $M$ to $M'$ 
extends as a meromorphic correspondence $F$ along any $CR$-path on $M$ 
and $F|_M[M]\subset M'$. If, moreover, $M'$ is strictly pseudoconvex 
then $F$ is a meromorphic map.
\end{thm}

Recall that a real hypersurface $M$ is called {\slsf minimal} if it does not
contain a non-constant germ of a complex hypersurface. 
The statement of this theorem is implicit  in \cite{SV} and is hidden inside of 
the proof of a stronger statement about {\slsf holomorphicity} of $F$. Unfortunately 
the proof of the holomorphicity of extension $F$ in \cite{SV} contains a gap. To make 
this point clear we give an outline of the proof of Theorem \ref{SVT}
referring step by step to \cite{SV}.

\proof The proof of the theorem breaks into several steps. Recall that a real submanifold
$\Sigma$ of $\cc^n$ is called {\slsf generic}  if its tangent space at any point contains
a complex subspace of minimal possible dimension. If $\dim_{\rr}\Sigma = 2n-2$ this means 
simply that $T_b\Sigma$ is not a complex subspace of $\cc^n$ for all $b\in \Sigma$. 
Genericity is obviously an open condition. 

\smallskip\noindent{\slsf Step 1.} The proof of Theorem \ref{SVT} can be reduced to the following
statement: {\it Let $\Omega$ be a domain in $M$ such that $f$ is meromorphic on $\Omega$ and let 
$b\in \partial \Omega$ be a boundary point such that $\Sigma \deff \partial \Omega$ is a
generic submanifold in a neighborhood of $p$. Then $f$ meromorphically extends to a neighborhood
of $b$.}

The proof is given in the section 4.1 of \cite{SV} and relies on the construction of special
ellipsoids from \cite{MP}. It is true that in section 4.1 of \cite{SV} $f$ is already supposed 
to be holomorphic on $\Omega$, but the proof goes through for any analytic objects, e.g., for 
meromorphic correspondences.

\smallskip Let $Q_b$ be the Segre variety of $M$ through $b$.   By Proposition 5.1 from \cite{S}
there exists a dense open subset $\omega$ of $Q_b$ such that for every $a\in \omega$ the
intersection $Q_a\cap \Omega$ is non-empty. Moreover, since $I_f$ is of codimension $\ge 2$
for a generic $a$ this intersection will be not contained in $I_f$. I.e., we can find 
$\xi \in \Omega \cap Q_a$ such that $f$ is holomorphic in a neighborhood $V_{\xi}$ of 
$\xi$. Let $V$ be a neighborhood of $Q_{\xi}$. 

\smallskip\noindent{\slsf Step 2.} {\it For an appropriate choice of $V_{\xi}$ and $V$ the set 
\begin{equation}
\eqqno(set-a-1)
A\deff \{(z,z')\in V\times \cc^{n'}: f(Q_z\cap V_{\xi})\subset Q'_{z'}\}
\end{equation}
is analytic in $V\times \cc^{n'}$, extends to an analytic set in $V\times \pp^{n'}$,
and this extension contains the graph of $f$ over $V_{\xi}$.}

Indeed, if $Q'$ is the defining polynomial of $M'$ the condition $f(Q_z\cap U_{\xi})\subset Q'_{z'}$
can be expressed as 
\begin{equation}
\eqqno(set-a-2)
Q'(f('t, h('t, \bar z), \bar z') =0.
\end{equation}
Here $'t = (t_1,...,t_{n-1})$ is a coordinate on the tangent plane $T_pM$ and $t_n=h('t, \bar z)$ 
is the equation of the Segre variety $Q_z$. After conjugation it is clear that this equation is 
holomorphic in $(z,z')$, see \cite{S,SV} for more details. Moreover, since \eqqref(set-a-2) is 
polynomial in $z'$ our set $A$ closes to an analytic set in $V\times \pp^{n'}$. This closure 
will be still denoted as $A$. The fact that $A$ contains the graph of $f$ over $V_{\xi}$
immediately follows from the invariance of Segre varieties: $f(Q_z)\subset Q'_{f(z)}$ whenever
everything is well defined (\ie sufficiently localized).

\smallskip Remark that dimension of $A$ might be bigger than $n$ simply because the set of 
$z'$ such that $Q'_{z'}\supset f(Q_z)$ can have a positive dimension for a fixed $z$. In order 
to reduce the size of $A$ proceed as follows. Consider the natural projections $\pi :A\to V$ 
and $\pi':A\to \pp^{n'}$. Remark that $\pi$ is proper simply because $\pp^{n'}$ is compact.
For an appropriate neighborhoods $V^*$ of $Q_a$ and $V_a$ of $a$ set
\begin{equation}
\eqqno(set-a*)
A^* \deff \{(z,z')\in V^*\times \pp^{n'}: \pi^{-1}(Q_z\cap V_a)\subset \pi'^{-1}(Q'_{z'})\}.
\end{equation}

\smallskip\noindent{\slsf Step 3.} {\it The set $A^*$ is analytic in $V^*\times \pp^{n'}$ of 
dimension $n$ and contains the graph of $f$ over (a shrinked, if necessary) $V_{\xi}$.}

The rough reason for $A^*$ to be of dimension $n$ is that both $\pi^{-1}(Q_z\cap V_a)$ and 
$\pi'^{-1}(Q'_{z'})$ are hypersurfaces in $A$ and therefore the set of $z'\in \pp^{n'}$ such that 
$\pi^{-1}(Q_z\cap V_a)\subset \pi'^{-1}(Q'_{z'})$ for a given $z\in V^*$ is generically finite. 
We refer to \S 3.2 of \cite{SV} for more details.

\smallskip As the graph of the correspondence $F$ which extends $f$ we take the irreducible component
of $A^*$ which contains the germ of the graph of $f$ over $\xi$. At points $z\in M\setminus I_F$
all values of $F$ are contained in $M'$ by unique continuation property of analytic functions
and by the connectivity of the graph of $F$.

\smallskip\noindent{\slsf Step 4.}{ \it If $M'$ is strictly pseudoconvex then $F$ is generically 
singlevalued, \ie is a meromorphic map.} 

\smallskip Suppose $z'\in F(z)$ for some $z\in M\setminus (I_F\cup R_F)$, where $R_F$ is the divisor 
of ramification of $F$. By the invariance of Segre varieties $F(Q_z)\subset Q_{z'}'$. This means that
for all branches $F_1,...,F_d$ of $F$ near $z$ one has $F_j(Q_z)\subset Q_{z'}'$. Let $z'=F_1(z)$ 
for simplicity and suppose that there is $w'=F_2(z)$ different from $z'$. Then $F_1(Q_z)\subset Q_{w'}'$
as well. But due to the strict pseudoconvexity of $M'$ we have that $Q_{z'}'\cap M' = \{z'\}$ 
and $Q_{w'}'\cap M'=\{w'\}$. I.e., the germs of $Q_{z'}'$ and $Q_{w'}'$ are disjoint. Contradiction, 
\ie $z'=w'$ and $F$ is singlevalued on $M\setminus (I_F\cup R_F)$. This implies that $R_F=\emptyset$ and $F$
is singlevalued on $M\setminus I_F$.

\smallskip\qed

\ifx\undefined\bysame
\newcommand{\bysame}{\leavevmode\hbox to3em{\hrulefill}\,}
\fi

\def\entry#1#2#3#4\par{\bibitem[#1]{#1}
{\textsc{#2 }}{\sl{#3} }#4\par\vskip2pt}

\end{document}